RATIONAL CUBOID. SEARCHING FOR A SOLUTION

B.S. Safin

February 2024

**Abstract**

This article discusses two versions of elliptic equations obtained from a system of equations describing a rational cuboid. Analysis of elliptic equations shows that they are equivalent, and that there are rational points on the elliptic curves, but they do not belong to any right triangle with rational sides.

Keywords: rational cuboid, elliptic equations, basic trigonometric identity for hyperbolic functions, congruent numbers.

**Preliminary remarks**

A rational cuboid (or a cuboid with all integer edges and diagonals) can be described by the following three quadratic equations:

$a^2 + b^2 = c^2$ $(x_1, y_1)$, $(x_{11}, y_{11})$ (1)

$a^2 + d^2 = e^2$ $(x_2, y_2)$, $(x_{21}, y_{21})$ (2)

$b^2 + d^2 = f^2$ $(x_3, y_3)$, $(x_{31}, y_{31})$ (3)

where **a, b, d** are the edges (let **a** and **b** be even), then the spatial diagonal is

$\mathbf{g^2 = a^2 + b^2 + d^2}$ (4)

where **c, e, f** are side diagonals. Obviously, the spatial diagonal $\mathbf{g^2}$ must be odd.

In two pairs of numbers $(x_1, y_1)$ and $(x_{11}, y_{11})$ each of the pairs are the numbers forming this right triangle (let c and b be odd), then it follows that:

$$x_1 = \sqrt{\frac{c+b}{2}},\ y_1 = \sqrt{\frac{c-b}{2}},\ x_{11} = \sqrt{\frac{c+a}{2}},\ y_{11} = \sqrt{\frac{c-a}{2}} \qquad x_{11} = \frac{(x_1+y_1)\cdot}{\sqrt{2}},\ y_{11} = \frac{(x_1 - y_1)}{\sqrt{2}} \qquad (4a)$$

and so on.

There is no common rule for determining the pairs $(x_1, y_1)$, $(x_{11}, y_{11})$, therefore in case of primitive triplets (a, b, c) we will take $(x_1, y_1)$ for integers. In case of non-primitive triplets specific equations and definitions shall be considered separately.

**Proposition 1**. For any right triangle $a^2 + b^2 = c^2$ $(x_1, y_1)$, $(x_{11}, y_{11})$ the following is true:

$$x^4_1 + y^4_1 + x^4_{11} + y^4_{11} = \left(\frac{3}{2}\right) \cdot c^2.$$

Proof: Raising each term to the fourth power, adding and combining similar terms we will arrive at the desired result.

A necessary condition for the existence of a cuboid is the presence of two even edges. This follows from the well-known fact that the sum of the squares of two odd numbers cannot be a square. Also, one of the three equations describing the cuboid must be divisible by 3, this is also a consequence of the long-proven fact that one of the legs of the Pythagorean equation must be divisible by 3.

Thus, in the system of equations (1-3) described above we will assume that **a, b, c** are even, therefore this triple of numbers cannot be primitive, just like the triple of numbers **b, d, f** of equation (3), it must be divisible by 3.

**We also assume that any cuboid can be reduced to this system (1-3).**

Consider the basic equation of the rational cuboid (4) $\mathbf{g^2 = a^2 + b^2 + d^2}$. This equation breaks down into three Pythagorean equations $\mathbf{b^2 + e^2 = c^2 + d^2 = a^2 + f^2 = g^2}$.

**Four representations of the spatial diagonal $\mathbf{g^2}$**

Each edge of the cuboid can be represented in four ways, which follows from equations (1-3) and relations (4a).

$\mathbf{a^2} = 4x_1^2y_1^2 = (x^2_{11} - y^2_{11})^2 = 4x_2^2y_2^2 = (x^2_{21} - y^2_{21})^2$

$\mathbf{b^2} = 4x^2_{11}y^2_{11} = (x^2_1 - y^2_1)^2 = 4x^2_3y^2_3 = (x^2_{31} - y^2_{31})^2$

$\mathbf{d^2} = 4x^2_{31}y^2_{31} = (x^2_3 - y^2_3)^2 = 4x^2_{21}y^2_{21} = (x_2^2 - y_2^2)^2$

**First case**

Note: since we do not know which two pairs of numbers $(x_1, y_1)$, $(x_{11}, y_{11})$ and $(x_2, y_2)$, $(x_{21}, y_{21})$ exactly form triangles (1) and (2), then we have to consider all four options.

$(x_1, y_1)$, and $(x_2, y_2)$; $(x_1, y_1)$ and $(x_{21}, y_{21})$; $(x_{11}, y_{11})$ and $(x_2, y_2)$; $(x_{11}, y_{11})$ and $(x_{21}, y_{21})$.

1) $a^2 + b^2 = c^2$ $(x_1, y_1)$, $(x_{11}, y_{11})$ – are the pairs of numbers forming this triangle

2) $a^2 + d^2 = e^2$ $(x_2, y_2)$, $(x_{21}, y_{21})$

3) $b^2 + d^2 = f^2$ $(x_3, y_3)$, $(x_{31}, y_{31})$

After checking all four options, we come to a single pair of numbers $(x_1, y_1)$ and $(x_2, y_2)$, which correspond to the basic equation of the cuboid $\mathbf{g^2 = a^2 + b^2 + d^2}$ in this configuration.

Since $a^2 = (4x_1^2y_1^2 = 4x_2^2y_2^2) \rightarrow +\mathbf{4x_1^2y_1^2 + (x^2_1 - y_1^2)^2 + (x_2^2 - y_2^2)^2} = g^2 = +\mathbf{4x_1^2y_1^2} + x_1^4 - \mathbf{2x_1^2y_1^2} + y_1^4 + x_2^4 -$
$- \mathbf{2x_2^2y_2^2} + y_2^4 = \mathbf{x_1^4 + y_1^4 + x_2^4 + y_2^4 = g^2} \rightarrow$ (5)

or, if we add and subtract $\mathbf{\pm 2x_1^2y_1^2} \rightarrow (x^2_1 + y_1^2)^2 + (x_2^2 - y_2^2)^2 = (x^2_1 - y_1^2)^2 + (x_2^2 + y_2^2)^2 =$
$g^2 = c^2 + d^2 = b^2 + e^2$.

**Second case**

Here, too, there are four options to consider.

$(x_1, y_1)$ and $(x_{11}, y_{11})$; $(x_1, y_1)$ and $(x_3, y_3)$; $(x_{11}, y_{11})$ and $(x_3, y_3)$ $(x_{11}, y_{11})$ and $(x_{31}, y_{31})$

1) $b^2 + a^2 = c^2$ $(x_1, y_1)$, $(x_{11}, y_{11})$ – are the pairs of numbers forming this triangle

2) $b^2 + d^2 = f^2$ $(x_3, y_3)$, $(x_{31}, y_{31})$

3) $a^2 + d^2 = e^2$ $(x_2, y_2)$, $(x_{21}, y_{21})$

After checking four pairs of numbers it appears, that to the basic equation of the cuboid $g^2 = a^2 + b^2 + d^2$ in this set corresponds the pair of numbers $(x_{11}, y_{11})$ and $(x_3, y_3)$, then the edges are $\mathbf{a^2} = (x^2_{11} - y^2_{11})^2$, $\mathbf{b^2} = 4x^2_{11}y^2_{11}$, $\mathbf{d^2} = (x^2_3 - y^2_3)^2 \rightarrow$

substitute $2x^2_{11}y^2_{11} = 2x^2_3y^2_3$ since $4x^2_{11}y^2_{11} = 4x^2_3y^2_3 \rightarrow$

$a^2 + b^2 + d^2 = x^4_{11} - \mathbf{2x^2_3y^2_3} + y^4_{11} + \mathbf{4x^2_{11}y^2_{11}} + x^4_3 - \mathbf{2x^2_3y^2_3} + y^4_3 \rightarrow \mathbf{x^4_{11} + y^4_{11} + x^4_3 + y^4_3 = g^2}$ (6)

**Third case**

Here, too, there are four options to consider.

$(x_3, y_3)$ and $(x_{31}, y_{31})$; $(x_3, y_3)$ and $(x_2, y_2)$ $(x_{31}, y_{31})$ and $(x_2, y_2)$ $(x_{31}, y_{31})$ and $(x_{21}, y_{21})$

1) $d^2 + b^2 = f^2$ $(x_3, y_3)$, $(x_{31}, y_{31})$ – are the pairs of numbers forming this triangle

2) $d^2 + a^2 = e^2$ $(x_2, y_2)$, $(x_{21}, y_{21})$

3) $b^2 + a^2 = c^2$ $(x_1, y_1)$, $(x_{11}, y_{11})$

In this case, the main equation corresponds in this set to the pair of numbers $(x_{21}, y_{21})$ and $(x_{31}, y_{31})$, and the edges are $\mathbf{a^2} = (x^2_{21} - y^2_{21})^2$, $\mathbf{b^2} = (x^2_{31} - y^2_{31})^2$, $\mathbf{d^2} = 4x^2_{31}y^2_{31} \rightarrow$

$g^2 = a^2 + b^2 + d^2 \rightarrow g^2 = (x^2_{21} - y^2_{21})^2 + (x^2_{31} - y^2_{31})^2 + 4x^2_{31}y^2_{31} \rightarrow$ since $(4x^2_{31}y^2_{31} = 2x^2_{31}y^2_{31} + 2x^2_{21}y^2_{21}) \rightarrow$

$x^4_{31} - \mathbf{2x^2_{31}y^2_{31}} + y^4_{31} + x^4_{21} - \mathbf{2\,x^2_{21}y^2_{21}} + y^4_{21} + \mathbf{4x^2_{31}y^2_{31}} \rightarrow$

$\rightarrow \mathbf{x^4_{21} + y^4_{21} + x^4_{31} + y^4_{31} = g^2}$ (7)

**Fourth case**

This case requires some clarification. Obviously, the solution to a Pythagorean triangle is the pairs of numbers that form this triangle, where each pair is essentially the coordinates of the intersection point of two hyperbolas in the Cartesian coordinate system (we are considering only the first quadrant for now). Therefore, each edge of a cuboid can be represented as a hyperbola. The solution to any Pythagorean triangle will be the intersection of two hyperbolas $(x \cdot y, x^2 - y^2)$ and if we rotate this pair of hyperbolas by 90 degrees (one clockwise and the other counterclockwise) we will get the second point of intersection with the second pair of numbers.

Example: $4^2 + 3^2 = 5^2$

The first pair of hyperbolas $(x \cdot y = 2, x^2 - y^2 = 3)$ $x = 2, y = 1$,

the second pair $(x \cdot y = 4, x^2 - y^2 = 1.5)$ $x = 3/\sqrt{2}, y = 1/\sqrt{2}$.

Now let us consider the equation numbered (3) $\mathbf{b^2 + d^2 = f^2}$ $(x_3, y_3)$, $(x_{31}, y_{31})$. One of the solutions to this equation is the intersection point of two hyperbolas $(x_3 \cdot y_3)$ and $(x^2_3 - y^2_3)$.

Let's draw the third hyperbola $a = (x^2_{11} - y^2_{11})$. This hyperbola will intersect the hyperbola $\mathbf{b = x_3 y_3}$ at the point $x_{11} = \frac{x_1 + y_1}{\sqrt{2}}$, $y_{11} = \frac{x_1 - y_1}{\sqrt{2}}$. Thus, on the hyperbola $\mathbf{b = x_3 y_3}$ there are two points with coordinates $(x_{11}, y_{11})$ and $(x_3, y_3)$.

Now we can use the properties of hyperbolic rotation and express the coordinates $(x_3, y_3)$ in terms of the coordinates $(x_{11}, y_{11})$.

$$x_3 = \frac{(x_1+y_1)\cdot\mathbf{m}}{\sqrt{2}}, \quad y_3 = \frac{(x_1 - y_1)}{\mathbf{m}\sqrt{2}}$$

where **m** is the hyperbolic rotation coefficient.

Let us express the edge "**a**" through the coordinates $x_3$, $y_3$. Since $x_3=(x_1+y_1)\cdot m/\sqrt{2}$, $y_3=(x_1-y_1)/m\cdot\sqrt{2}$.

Then $(x_1+y_1) = 2x_3/m\cdot\sqrt{2}$, $(x_1-y_1) = 2y_3\cdot m/\sqrt{2}$

Adding and subtracting these expressions we get:

$x_1=(x_3+y_3m^2)/m\cdot\sqrt{2}$, $y_1=(x_3-y_3m^2)/m\cdot\sqrt{2}$, and $\mathbf{a}=2x_1y_1 = (x_3^2-y_3^2m^4)/m^2$.

We also take $x_3/y_3=\text{ch}\alpha \rightarrow$ let $\mathbf{a_1}=(x_1-y_1)/(x_1+y_1) \rightarrow \mathbf{m^2=a_1\cdot ch\alpha}$. (8)

Now we can express the spatial diagonal as: $\mathbf{g^2=f^2+a^2}=(x_3^2+y_3^2)^2+(x_3^2-y_3^2m^4)^2/m^4$,

By opening the brackets and bringing to a common denominator, we obtain

$$\mathbf{g^2=\frac{(m^4+1)(x_3^4+y_3^4m^4)}{m^4}} \quad (9)$$

For each of the three sets we received a representation of the spatial diagonal $\mathbf{g}^2$.

$$x_1^4+y_1^4+x_2^4+y_2^4=g^2 \quad (5)$$

$$x^4_{11}+y^4_{11}+x^4_3+y^4_3=g^2 \quad (6)$$

$$x^4_{21}+y^4_{21}+x^4_{31}+y^4_{31}=g^2 \quad (7)$$

The fourth representation is a variant of the first one expressed in terms of $(x_3, y_3)$.

$$x^4_3 + \frac{x^4_3}{m^4} + m^4 y^4_3 + y^4_3 = (m^4+1)\cdot\frac{x_3^4+y_3^4m^4}{m^4} = g^2 \quad (9)$$

**First elliptic equation**

Let us consider the equation (5) $x_1^4+y_1^4+x_2^4+y_2^4=g^2$. This equation refers to the first case in which we found out that pairs of numbers $(x_1, y_1)$ and $(x_2, y_2)$ are suitable for this equation, that is,

$\mathbf{a}=2x_1y_1=2x_2y_2 \rightarrow \frac{x_2}{x_1} = \frac{y_1}{y_2} = \mathbf{t} \rightarrow \frac{x_2\mathbf{k}}{x_1\mathbf{k}} = t \rightarrow$

$y_1=x_2\mathbf{k}=x_1\mathbf{t\cdot k}$, $y_2= x_1\mathbf{k} \rightarrow x_2=x_1t$, $y_1=y_2\mathbf{t}=x_2\mathbf{k}=x_1\mathbf{kt}$, $x_1=x_1$ (we take $x_1<y_1$) (9a)

where **k** and **t** are also the hyperbolic rotation coefficients like **m.**

Next, if we take $y_1^4=\frac{x_2^4y_2^4}{x_1^4}$ and substitute it into Eq. (5), we obtain

$(x_1^4+x_2^4)\cdot(x_1^4+y_2^4) = g^2x_1^4$ (10)

and if we take $y_2^4=\frac{x_1^4y_1^4}{x_2^4}$ and substitute it into (5), then we obtain

$(x_1^4+x_2^4)\cdot(x_2^4+y_1^4) = g^2x_2^4$ (11)

Equations (10) and (11) are equivalent to each other, so we make use of take any of them.

From equation (11) follows that $(x_1^4+x_2^4)\cdot(x_2^4+y_1^4) = g^2x_2^4$, and taking into account that $x_2^4=x^4_1t^4$ we obtain:

$x_1^4t^4+y_1^4+x_1^4t^8+y_1^4t^4=g^2t^4$.

Dividing this equation by $t^2$, and taking into account that $\left(\frac{y_1^4}{t^2}=y_1^2y_2^2\right)$ we get

$t^2(x_1^4+y_1^4) + x_1^4t^6+y_1^2y_2^2=g^2t^2$.

By moving the first two terms to the right side, and taking into account that $(g^2-x_1^4-y_1^4=x_2^4+y_2^4)$, $x_2=x_1t$, $y_2=x_1k$, we obtain:

$y_1^2y_2^2 = t^2\cdot(g^2-x_1^4-y_1^4) - x_1^4t^6$.

Multiplying this expression by $x^2_1 \rightarrow x^2_1y_1^2y_2^2= x^2_1t^2(x_2^4+y_2^4) - x_1^6t^6$.

Let's take $\mathbf{y^2}=x^2_1y_1^2y_2^2$ and $\mathbf{x}=x^2_1\mathbf{t^2}$, consequently

$\mathbf{y^2=x\cdot(x_2^4+y_2^4)-x^3}$ or $\mathbf{y^2=x^3-(x_2^4-y_2^4)\cdot x}$. (12)

(the equations are equivalent, which can be checked by opening the brackets).

As we found out earlier, the main equation (4) breaks down into 3 equations:

$b^2+e^2 = c^2+d^2 = a^2+f^2 = g^2$.

For the cuboid to be perfect, it is necessary that all these three equations have rational or integer solutions. Knowing the known cuboids, one can easily construct the elliptic equations in explicit form. Let's consider the elliptic equations for the first case $\mathbf{y^2=x\cdot(x_2^4+y_2^4)-x^3}$, or $\mathbf{y^2=x^3-(x_2^4-y_2^4)\cdot x}$.

We believe that we can always reduce the system of equations (1-3) to the first case, therefore, without loss of generality, we can consider only this case.

For example, for a cuboid with edges (104, 153, 672) $x_1=2$, $y_1=26$ and $x_2=13$, $y_2=4$, $\mathbf{y}^2=x^2_1y_1^2y_2^2$, $\mathbf{x}=x_2^2$,

| | | |
|---|---|---|
| $y^2=x^3-28305\cdot x$ | $(x=13^2, y=2\cdot 26\cdot 4)$ | – here we took the difference for $(x_2^4-y_2^4)$. |
| $y^2=x\cdot 28817-x^3$ | $(x=4^2, y=26^2)$ | – here we swap $x_2$ and $y_2$. |
| $y^2=x\cdot 456992-x^3$ | $(x=2^2, y=2\cdot 26^2)$ | – here we swap the pair $(x_2, y_2)$ and $(x_1, y_1)$. |
| $y^2=x^3-456960\cdot x$ | $(x=26^2, y=4\cdot 26)$ | – here we took the difference $(x_1^4-y_1^4)$. |

Equations (6), (7), (9) were not considered due to their equivalence to the equation (5).

Considering that $\rightarrow \frac{x_2}{x_1}=\frac{y_1}{y_2}=\mathbf{t}\rightarrow \frac{x_2\mathbf{k}}{x_1\mathbf{k}}=t\rightarrow y_1=x_2\mathbf{k}=x_1\mathbf{t}\cdot\mathbf{k}$ and $x_1\mathbf{k}=y_2$, we can divide the equations $y^2=x^3-(x_2^4-y_2^4)\cdot x$ and $y^2=(x_2^4+y_2^4)\cdot x-x^3$ by $x^6_1$, so we get the identities:

$k^4t^2=t^2(k^4+t^4)-t^6$ or $k^4t^2=t^6+t^2(k^4-t^4)$

that are valid for any **t** and **k.** (here we take $x_1<y_1$, if $x_1>y_1$ then we will need to divide by $\mathbf{y_1^6}$).

The examples show that there are integer points on the elliptic curve (12), therefore, there are also rational ones. The question arises: is it possible to construct a right triangle with rational or integer sides based on these points?

**Lemma 1.** If the numbers $\mathbf{x_2^4}$, $\mathbf{y_2^4}$ are integers, then on elliptic curves of the form $\mathbf{y^2=x^3-(x_2^4-y_2^4)\cdot x}$ and $\mathbf{y^2=(x_2^4+y_2^4)\cdot x-x^3}$ there is no point from which a right triangle could be with integer or rational sides.

Proof: Let's use well-known formulas that allow constructing a right triangle with rational or integer sides based on a point on an elliptic curve (Koblitz **[1]**). Let a triangle be $a^2+b^2=c^2$, then:

$$a=\frac{x^2-n^2}{y},\ b=\frac{2nx}{y}\ \text{and}\ c=\frac{x^2+y^2}{y}$$

in our case $x=\mathbf{x_2^2}$, $y=\mathbf{x_1y_1y_2}$ and $n^2=\mathbf{x_2^4+y_2^4}$. Let's construct a triangle →

$$\frac{y^8_2}{x^2_1y_1^2y_2^2}+4\cdot\frac{x_2^4(x_2^4+y_2^4)}{x^2_1y_1^2y_2^2}=\frac{(2\,x_2^4+y_2^4)^2}{x^2_1y_1^2y_2^2}$$

Reducing both sides of the equation by $\mathbf{(x^2_1y_1^2y_2^2)}$ we get $\mathbf{y^8_2+4\cdot x_2^4(x_2^4+y_2^4)=(2\cdot x_2^4+y_2^4)^2}$ (13)

Obviously, one leg $\mathbf{y^8_2}$ and the hypotenuse $\mathbf{(2\cdot x_2^4+y_2^4)^2}$ are squares. Let's consider the second leg $\mathbf{4\cdot x_2^4(x_2^4+y_2^4)}$. Since in our case we are considering a system of equations 1-3, then the lateral hypotenuse **e** is an odd number, therefore, the numbers $\mathbf{x_2^4}$, $\mathbf{y_2^4}$ have different parities. Then we can use a well-known **theorem** from the number theory which states the following: «**If the product of two relatively prime natural numbers is the nth power, then each of the factors will also be the nth power**». In our particular case, the nth power is a square. **[7]**. $\mathbf{4\cdot x_2^4}$ is even, and $\mathbf{(x_2^4+y_2^4)}$ is odd. We can apply the theorem. Consequently it turns out that in order for the product of $\mathbf{4\cdot x_2^4(x_2^4+y_2^4)}$ to be a square, it is necessary that $\mathbf{(x_2^4+y_2^4)}$ also be a square, but Fermat's theorem does not allow this for the case of $n=4$. Therefore we cannot construct a right triangle if the numbers $\mathbf{x_2^4}$, $\mathbf{y_2^4}$ are integers.

**Second elliptic equation**

To obtain the first elliptic equation, we actually used only the numbers that form two of the three equations that describe the cuboid, that is $(x_1, y_1)$ and $(x_2, y_2)$. Those were the equations $a^2+b^2=c^2$ and $a^2+d^2=e^2$. Although indirectly the numbers $(x_3, y_3)$ were also present there, since we carried out transformations of the equation (4) $\mathbf{g^2=a^2+b^2+d^2}$, in which, if we take the terms on the right in pairs, we get three Pythagorean equations and one of them will be the equation $a^2+f^2=g^2$.

To construct the second elliptic equation, we will use the third equation of the system $b^2+d^2=f^2$ $(x_3, y_3)$, $(x_{31}, y_{31})$, as well as the equation $a^2+f^2=g^2$.

To obtain the second elliptic equation we need two parametrizations of the Pythagorean equation: let **x** and **y** be the numbers forming the Pythagorean triangle:

**Proposition 2.** $(\text{ch}\alpha=x/y)\rightarrow 4y^4\text{ch}^2\alpha+y^4\text{sh}^4\alpha=y^4(1+\text{ch}^2\alpha)^2)$ (14)

**Proposition 3.**

$$\mathbf{a_1}=\frac{x-y}{x+y}\rightarrow y^4\frac{4(1+a_1)^2}{(1-a_1)^2}+y^4\frac{16a^2_1}{(1-a_1)^4}=y^4\frac{4(1+a^2_1)^2}{(1-a_1)^4}\ \text{in our case}\ \mathbf{a_1}=\frac{x_1-y_1}{x_1+y_1}=\frac{x_1(1-kt)}{x_1(1+kt)}\quad(15)$$

$y^4$ we find from the relation $a^2=4x^2_1y^2_1=y^4\frac{4(1+a_1)^2}{(1-a_1)^2}$, since we have taken $y_1=x_1kt$ (see 9a), then we have $y^4=x^4_1$ (if we take $x_1=y_1kt$, then $y^4=y^4_1$). The parameterization of equations (14) and (15) need not be proved, since everything can be checked by direct substitution. In the first case, we come to the identity $sh^2\alpha = sh^2\alpha \rightarrow sh^2\alpha = ch^2\alpha - 1 \rightarrow ch^2\alpha - sh^2\alpha = 1$, and in the second case to the identity $(1+a_1^2) = (1+a_1^2)$. It is important to note that these parameterizations do not give us an idea of the triangle itself, that is, we cannot say whether the squares of the legs and hypotenuse are rational, integer or irrational.

Let us consider the equation $b^2+d^2=f^2$ $(x_3, y_3)$, $(x_{31}, y_{31})$.

Parameterization (14) brings this equation to the form $4y_3^4ch^2\alpha + y_3^4sh^4\alpha = y_3^4(1+ch^2\alpha)^2$.

Let us recall the hyperbolic rotation coefficient (8) $\mathbf{m^2=a_1 \cdot ch\alpha}$.

If the hyperbola $\mathbf{a=x_1y_1}$ is rotated by 90 degrees, then it will intersect the hyperbola $\mathbf{b=x_3y_3}$ at the point $(x_{11}, y_{11})$, then $\mathbf{m}=\frac{x_3}{x_{11}}=\frac{y_{11}}{y_3}$ is the hyperbolic rotation coefficient.

Next $\rightarrow x_3=x_{11}\cdot\mathbf{m}=(x_1+y_1)\cdot\frac{\mathbf{m}}{\sqrt{2}}$, $y_3=\frac{y_{11}}{\mathbf{m}}=\frac{x_1 - y_1}{\mathbf{m}\cdot\sqrt{2}}$.

Also let $ch\alpha=\frac{x_3}{y_3}$, $a_1=\frac{x_1 - y_1}{x_1+y_1}$, then $ch\alpha=\frac{x_3}{y_3}=\frac{m^2}{a_1}$ or $\mathbf{m^2 = a_1 \cdot ch\,\alpha}$.

Now ch α in a triangle $b^2+d^2=f^2=4y_3^4ch^2\alpha+y_3^4sh^4\alpha=y_3^4(1+ch^2\alpha)^2$ can be expressed as $\mathbf{ch\alpha=m^2/a_1}$.

Consider the equation $\mathbf{a^2+f^2=g^2}$.

Parameterization (15) gives $a^2=x^4_1\frac{4(1+a_1)^2}{(1-a_1)^2}$, parameterization (14) gives $f^2=y^4_3(1+ch^2\alpha)^2$.

To express the spatial diagonal $g^2$ we take relation (9), where $g^2=(m^4+1)\cdot\frac{x_3^4+y_3^4m^4}{m^4} \rightarrow \mathbf{a^2+f^2=g^2} \rightarrow$

$$x^4_1\frac{4(1+a_1)^2}{(1-a_1)^2} + y^4_3(1+ch^2\alpha)^2=(m^4+1)\cdot\frac{x_3^4+y_3^4m^4}{m^4}.$$

Let us transform $\mathbf{y_3}=\frac{(x_1-y_1)}{m\sqrt{2}} \rightarrow \left[y_3^4 = x^4_1\cdot\frac{(1-kt)^4}{4\cdot m^4}\right]$

Here we have taken $x^4_1$ out of brackets taking into account that in our system of equations $y_1=x_1\cdot kt$.

Then the equation $\mathbf{a^2+f^2=g^2}$ takes the form:

$$x^4_1\frac{4(1+a_1)^2}{(1-a_1)^2} + x^4_1\cdot\frac{(1-kt)^4}{4\cdot m^4}\cdot(1+ch^2\alpha)^2=(m^4+1)\cdot\frac{x_3^4+y_3^4m^4}{m^4}.$$

Let us express $(x_1-y_1)^4= x^4_1(1-kt)^4$ in terms of $\mathbf{a_1}$. Since $\rightarrow y_1=x_1\mathbf{t\cdot k}\rightarrow \mathbf{a_1}=\frac{x_1-y_1}{x_1+y_1}=\frac{x_1(1-kt)}{x_1(1+kt)}$, we reduce them by $x_1$. Now kt may be expressed as $kt=\frac{1+a_1}{1-a_1}$, a $\rightarrow (1-kt)^4=\left\{1-\frac{(1+a_1)^4}{(1-a_1)^4}\right\}=16\cdot\frac{a^4_1}{(1-a_1)^4}$. So, the equation $\mathbf{a^2+f^2=g^2}$ takes the following form:

$$\mathbf{x^4_1\frac{4(1+a_1)^2}{(1-a_1)^2}+x^4_1\cdot\frac{16\,a^4_1}{4\cdot m^4(1-a_1)^4}(1+ch^2\alpha)^2=(m^4+1)\cdot\frac{x_3^4+y_3^4m^4}{m^4}}.$$

We also transform the expression (9)$\rightarrow (m^4+1)\cdot\frac{x_3^4+y_3^4m^4}{m^4}=g^2$.

Let's get rid of the coordinates $(x_3^4, y_3^4)$.

$$ch^4\alpha=\frac{x_3^4}{y_3^4}\rightarrow(m^4+1)\cdot\frac{y_3^4\cdot ch^4\alpha+y_3^4m^4}{m^4}\rightarrow g^2=y_3^4(m^4+1)\cdot\frac{ch^4\alpha+m^4}{m^4}$$

$\left[y_3^4 = x^4_1\cdot\frac{(1-kt)^4}{4\cdot m^4}\right]$ since $(1-kt)^4=16\cdot\frac{a^4_1}{(1-a_1)^4}$ we substitute these values of $y_3^4$ and $(1-kt)^4$ into the expression:

$$g^2=y_3^4(m^4+1)\cdot\frac{ch^4\alpha+m^4}{m^4}=\frac{g^2}{x^4_1}=(m^4+1)\cdot16\frac{a_1^4(ch^4\alpha+m^4)}{4m^8(1-a_1)^4}\rightarrow 4a_1^4(a_1^2ch^2\alpha+1)\cdot\frac{ch^2\alpha(ch^2\alpha+a^2_1)}{a^4_1ch^4\alpha(1-a_1)^4}$$

Here we replaced $m^4$ with $a_1^2\cdot ch^2\alpha$, $m^8$ with $a_1^4\cdot ch^4\alpha$ and reduced by $(4a_1^4)$ and by $ch^2\alpha \rightarrow$

$$\frac{g^2}{x^4_1}=\frac{4(a_1^2ch^2\alpha+1)(a_1^2+ch^2\alpha)}{ch^2\alpha(1-a_1)^4}$$

After replacing chα back with $m^2/a_1$ and minor transformations, we get:

$$\frac{g^2}{x_1^4}=4\cdot\frac{(m^4+1)(a^4_1+m^4)}{m^4(1-a_1)^4}.$$

Now the equation $a^2+f^2=g^2$ has taken its final form:

$$\mathbf{x^4_1\frac{4(1+a_1)^2}{(1-a_1)^2}+x^4_1\cdot16\cdot\frac{a^4_1}{4\cdot m^4(1-a_1)^4}\cdot(1+ch^2\alpha)^2=4\cdot x^4_1\frac{(m^4+1)(a^4_1+m^4)}{m^4(1-a_1)^4}}$$

We reduce this equation by $4x^4_1$ and solve it for $ch^2\alpha$. →
$(1+a_1)^2\cdot(1-a_1)^2\cdot m^4+(1+ch^2\alpha)^2\cdot a^4_1=(m^4+1)(a^4_1+m^4)\to$
$m^4-2a^2_1\cdot m^4+a^4_1\cdot m^4+a^4_1\cdot(1+ch^2\alpha)^2=a^4_1\cdot m^4+m^8+a^4_1+m^4$.
After combining the similar terms and reducing we get → $a^4_1(1+ch^2\alpha)^2=(m^4+a^2_1)^2$, since according to (8) $m^4=a^2_1\cdot ch^2\alpha$,
then $a^4_1(1+ch^2\alpha)^2=(a^2_1\cdot ch^2\alpha+a^2_1)^2\to 1+ch^2\alpha=1+ch^2\alpha\to ch^2\alpha=ch^2\alpha=1+sh^2\alpha\to$ $\mathbf{ch^2\alpha-sh^2\alpha=1}$.
From this identity we obtain the second elliptic equation.
Since $ch^2\alpha=\frac{m^4}{a^2_1}\to\frac{m^4}{a^2_1}-1=sh^2\alpha$.
Multiply both sides by $m^2\to\frac{m^6}{a^2_1}-m^2=m^2\cdot sh^2\alpha\to a^2_1\cdot m^2=m^6-a^2_1m^2\cdot sh^2\alpha$.
($\mathbf{m^2}=a_1\cdot ch\alpha\to\mathbf{m^6}=a^3_1ch^3\alpha$).
$\mathbf{a^2_1\cdot m^2=a^3_1ch^3\alpha-a^2_1sh^2\alpha\cdot a_1ch\alpha}\to$, ($y=a_1\cdot m$, $x=a_1ch\alpha$, $n=a_1sh\alpha$) $\mathbf{y^2=x^3-n^2x}$ (16)

**Example**: Let's take a cuboid with edges (44, 117, 240). Let us bring it to the form that describes **the first case.** We believe that we can always bring the system to this case.
Then the following relations must be true: $\mathbf{a}^2=4x_1^2y_1^2=4x_2^2y_2^2$, and $x_1^4+y_1^4+x_2^4+y_2^4=g^2$
That's why **a**= 240, **b**=44, **d**=117. The system of equations describing this cuboid will look like follows:
$240^2+44^2=244^2\quad(x_1=12, y_1=10, x_{11}=11\sqrt{2}, y_{11}=\sqrt{2})$
$240^2+117^2=267^2\ (x_2=8\sqrt{3}, y_2=5\sqrt{3}, x_{21}=\frac{13\sqrt{3}}{\sqrt{2}}, y_{21}=\frac{3\sqrt{3}}{\sqrt{2}})$
$44^2+117^2=125^2\ \left(x_3=11, y_3=2, x_{31}=\frac{13}{\sqrt{2}}, y_{31}=\frac{9}{\sqrt{2}}\right)$
Let's calculate what $a_1$, m, $ch\alpha$, $sh\alpha$ are equal to.
$a_1=\frac{x_1-y_1}{x_1+y_1}=\frac{12-10}{12+10}=\frac{1}{11}$
$ch\alpha=\frac{x_3}{y_3}=\frac{11}{2}=5{,}5$
then $sh^2\alpha=29.25$, $m^2=a_1\cdot ch\alpha=0.5$.
Now equation (16) will have the following form:
$a^2_1\cdot m^2=a^3_1ch^3\alpha-a^2_1\cdot sh^2\alpha\cdot a_1ch\alpha\to\mathbf{y^2=x^3-14157x}$,
where ($\mathbf{x}=11^2$, $\mathbf{y}=242$).
For the first elliptic equation (12) $\mathbf{y^2=x\cdot(x_2^4+y_2^4)-x^3}$ we get $y^2=42489\cdot x-x^3$ (where $y^2=108000$, $x=192$). For the first case you can also obtain other equations, as it was shown above in the example for the first elliptic equation.
For our equation (16) $\mathbf{n^2}=a^2_1\cdot\sin^2\alpha$, and we cannot say whether this number can be congruent or not, just like about the numbers **x** and **y.** Here we need another criterion to determine the possibility of constructing a right triangle with rational sides.
Let us construct a right triangle from the resulting equation (16). We use the formulas with which we built a triangle from the first elliptic equation:
$a=\frac{x^2-n^2}{y}$, $b=2n\cdot\frac{x}{y}$ and $c=\frac{x^2+y^2}{y}$
In our case $\mathbf{x}=a_1\cdot ch\alpha$, $\mathbf{y}=a_1\cdot m$, $\mathbf{n}=a_1\cdot sh\alpha$, so we get:

$$\frac{a^2_1}{m^2}+\frac{4\,a^2_1\cdot sh^2\alpha\cdot ch^2\alpha}{m^2}=\frac{a^2_1(sh^2\alpha+ch^2\alpha)^2}{m^2}$$

Reducing both sides by $\frac{a^2_1}{m^2}$, we obtain $\mathbf{1+4sh^2\alpha\cdot ch^2\alpha=(sh^2\alpha+ch^2\alpha)^2}\to ch^2\alpha-sh^2\alpha=1$ (17)
**Proposition 4.** If the numbers $sh^2\alpha$ and $ch^2\alpha$ are integers, then using this equation (16) we will not be able to construct a right triangle.
Proof: The integers $sh^2\alpha$ and $ch^2\alpha$ must have different parities by definition.
Since $ch^2\alpha-sh^2\alpha=1$, but the difference of the squares of two integers cannot be equal to one, therefore, one of the numbers $ch^2\alpha$ or $sh^2\alpha$ is not a square, then, according to **the theorem** from the first case, the term $4sh^2\alpha\cdot ch^2\alpha$ also cannot be a square.

**Equivalence of two equations**

Let us prove that equation (12) $\mathbf{y^2=x^3-(x_2^4-y_2^4)\cdot x}$ and equation (16) $\mathbf{a^2{}_1\cdot m^2=m^6-a^2{}_1m^2sh^2\alpha}$ are equivalent to each other. To do this, we will show that from the basic equation describing the cuboid $\mathbf{g^2=a^2+b^2+d^2}$ we can obtain not only the equation (12) $\mathbf{y^2=\pm x^3-(x_2^4\pm y_2^4)\cdot x}$ but also the equation (15) $\mathbf{a^2{}_1\cdot m^2=a^3{}_1ch^3\alpha-a^2{}_1\cdot sh^2\alpha\cdot a_1ch\alpha=}$ (where $\mathbf{x}=a_1\cdot ch\alpha$, $\mathbf{y}=a_1\cdot m$, $\mathbf{n}=a_1\cdot sh\alpha$).

**Lemma 2**. Equation (12) $\mathbf{y^2=\pm x^3-(x_2^4\pm y_2^4)\cdot x}$ and equation (16) $\mathbf{a^2{}_1\cdot m^2=m^6-a^2{}_1m^2sh^2\alpha}$ are equivalent to each other and are equivalent to the main trigonometric identity for hyperbolic functions

$\mathbf{ch^2\alpha - sh^2\alpha=1.}$

Proof:

Parameterization (15) gives us: $\mathbf{a^2}=4\,x^4{}_1\cdot\frac{(1+a_1)^2}{(1-a_1)^2}$, и $\mathbf{b^2}=x^4{}_1\cdot\frac{16a^2{}_1}{(1-a_1)^4}$,

Parameterization (14) gives us: $\mathbf{d^2}=y^4{}_3\cdot sh^4\alpha=x^4{}_1\cdot 16\cdot\frac{a^4{}_1}{4\cdot m^4(1-a_1)^4}\cdot sh^4\alpha$,

$\mathbf{g^2}=x^4{}_1\frac{(m^4+1)(a^4{}_1+m^4)}{m^4(1-a_1)^4}$ - this expression was obtained earlier.

Now the equation $a^2+b^2+d^2=g^2$ can be rewritten as follows::

$$4\,x^4{}_1\cdot\frac{(1+a_1)^2}{(1-a_1)^2}+x^4{}_1\cdot\frac{16a^2{}_1}{(1-a_1)^4}+16\,x^4{}_1\cdot a^4{}_1\cdot\frac{sh^4\alpha}{4m^4(1-a_1)^4}=4x^4{}_1\cdot\frac{(m^4+1)(m^4+a_1^4)}{m^4\cdot(1-a_1)^4} \quad (18)$$

Reducing this equation by $4\cdot x^4{}_1$ and solving it for $sh^4\alpha$, we get

$$sh^4\alpha=\left(\frac{m^4-a_1{}^2}{a_1{}^2}\right)^2=sh^4\alpha$$

since $ch^2\alpha=\frac{m^4}{a_1{}^2}$.→ $sh^2\alpha=sh^2\alpha\rightarrow ch^2\alpha-sh^2\alpha=1$.

We have already obtained from this identity the elliptic equation (16) $\mathbf{y^2=x^3-n^2x}$. Now we can say that the equivalence of the pair of elliptic equations $\mathbf{y^2=x^3-(x_2^4-y_2^4)\cdot x}$ and $\mathbf{y^2=x\cdot(x_2^4+y_2^4)-x^3}$ with ($\mathbf{y}=x_1y_1y_2$, a, $\mathbf{x}=x^2{}_1t^2$, $\mathbf{n^2}=x_2^4\pm y_2^4$) to the equation $\mathbf{y^2=x^3-n^2x}$ with ($\mathbf{x}=a_1\cdot ch\alpha$, $\mathbf{y}=a_1\cdot m$, $\mathbf{n^2}=a^2{}_1\cdot sh^2\alpha$) has been proven.

**Discussion of the results**

Since we believe that in the system of equations (1-3) all numbers are integers, and that the equations describing the cuboid can always be reduced to this system, then the first triple of numbers in equation (1) **a, b, c** is not a primitive triple, since the numbers (**a, b**) are even by definition. Also, the triple of numbers from equation (3) (**b, d, f**) also cannot be primitive by definition, since it must be divisible by the number 3. Thus, only the numbers from the second equation (**a, d, e**) may possibly constitute a primitive triple. If this triple is primitive, then one pair of numbers forming it ($x_2, y_2$) or ($x_{21}, y_{21}$) must be integers. Then it turns out that according to equation (13) $\mathbf{y^8{}_2+4\cdot x_2^4(x_2^4+y_2^4)=(2\cdot x_2^4+y_2^4)^2}$ we cannot construct a right triangle with integer or rational sides, due to the fact that the expression $\mathbf{4\cdot x_2^4(x_2^4+y_2^4)}$ cannot be a square according to the theorem stipulating that "**if the product of two relatively prime natural numbers is the nth power, then each of the factors will also be the nth power**." In our case n=2. [7], and $(x_2^4+y_2^4)$ cannot be a square due to Fermat's theorem for power 4, hence the leg $4\cdot x_2^4(x_2^4+y_2^4)$ is not a square too. If the triple is not primitive, then since our hypotenuse **"e"** by condition must be an integer and odd number, then the numbers forming the triangle $a^2+d^2=e^2$ ($x_2, y_2$), ($x_{21}, y_{21}$) must have different parities, and even if they are not integers, then the squares of these numbers must be integers, since all solutions of the Pythagorean equation in terms of integers are determined by the formulas: $\mathbf{a}=2\cdot x\,y\cdot p$, $\mathbf{d}=(x^2-y^2)\cdot p$, $\mathbf{e}=(x^2+y^2)\cdot p$ (we assume that **a** is even, **d, e** are odd), where x, y, p are natural numbers [3]. By definition $x_2=\sqrt{\frac{(e+d)\cdot p}{2}}$, $y_2=\sqrt{\frac{(e-d)\cdot p}{2}}$, then $x_2=x\cdot\sqrt{p}$ and $y_2=y\cdot\sqrt{p}$, thus the numbers $\mathbf{x^2{}_2}$ and $\mathbf{y^2{}_2}$ must be integers, since x, y, p are natural numbers (we assume that $x > y$).

The same reasoning can be applied to equations (1, 3) in our system, which we are considering, although the triplets of these equations are not primitive.

Therefore, if in equation (13) we substitute ($x_2, y_2$) for ($x_1, y_1$), the result will be the same.

It is not possible to construct a right triangle with integer or rational sides. Since we have shown that equations of the form (12) $\mathbf{y^2=x\cdot(x_2^4\pm y_2^4)\pm x^3}$ and equations of the form (16) $\mathbf{a^2{}_1\cdot m^2=a^3{}_1ch^3\alpha-a^2{}_1sh^2\alpha\cdot a_1ch\alpha}$ are equivalent to each other, then all arguments can be applied to both equations. Thus, even if there are integer points on the elliptic curves (12) and (16), it is not possible to construct a right triangle with integer or rational sides.

**Ways for obtaining a rational cuboid**

I know of two methods for obtaining cuboids with one or more non-rational elements. The first method was invented by Euler. He found formulas for obtaining a cuboid with a non-integer diagonal ($X=n^6-15n^4+15n^2-1, Y=6n^5-20n^3+6n, Z=8n^5-8n$) where $n \geq 2$ is a natural number **[2].** But in fact, these formulas do not describe all cuboids. Cuboids with one non-integer element are also called Euler bricks.

The idea of the second method belongs to Shedd (Ch. L. Schedd), and was supplemented by Serpinsky (Serpinsky V.) **[3].** Generalizing their ideas, the method can be called the multiplication of Pythagorean triangles. The essence of the method is that the Brahmagupta -Fibonacci identity is used to multiply the legs of two Pythagorean triangles. If the triangles are primitive, then at the output we will get two primitive triangles. Another pair of triangles can be obtained by multiplying each original triangle by the hypotenuse of the other triangle. By combining these triangles, you can get 4 cuboids with an integral spatial diagonal and 4 cuboids with a non-integer spatial diagonal. Some of the elements of a cuboid may be negative and non-integer. The result is a rather exotic algebraic surface. The next paragraph gives an example of 8 cuboids obtained by multiplying two primitive Pythagorean equations. More details about the method can be found in **[6].** How to find pairs of Pythagorean triangles, the multiplication of which produces a cuboid with one non-integer element, is an open question.

The publication [2] describes a way for obtaining cuboids with intact edges and a spatial diagonal. Also a lot of information concerning cuboids is contained in [3].

**An interesting fact about the rational cuboid**

**Proposition 5.** The fourth power of the spatial diagonal ($g^4$) of a cuboid is equal to the difference of the sum of the fourth powers of the side diagonals of the cuboid and the fourth powers of the edges, that is: $g^4=c^4+e^4+f^4-a^4-b^4-d^4$.

Proof: Let's square both sides of the equation $(g^2)^2=(a^2+b^2+d^2)^2$, after simple transformations we get the desired result.

**Lemma 3**: The fourth power of the spatial diagonal ($g^4$) is equal to the difference of the quadruple square of the area of the triangle built on the side diagonals, and the quadruple square of the area of the triangle built on the edges of the cuboid, that is $\mathbf{g^4=16 \cdot S_1^2-16 \cdot S_2^2}$.

Proof: g is the spatial diagonal, c, e, f are the side diagonals, a, b, d are the edges.

Using Heron's formula for the area of a triangle, we calculate the area of the triangle formed by the side diagonals c, f, e:

$$S_1=\left(\frac{1}{4}\right)\cdot\sqrt{(c^2+e^2+f^2)^2-2(c^4+e^4+f^4)}$$

but $c^2+e^2+f^2=2g^2$, and $c^4+e^4+f^4=g^4+a^4+b^4+d^4$, so the radical expression can be rewritten as follows:

$16 \cdot S_1^2=(2g^2)^2-2(c^4+e^4+f^4) \to 16 \cdot S_1^2=4 \cdot g^4-2(g^4+a^4+b^4+d^4) \to 2g^4-2(a^4+b^4+d^4)=g^4+[g^4-$
$-2(a^4+b^4+d^4)]$,

where $16 \cdot S_2^2=g^4-2(a^4+b^4+d^4)$ is the square of the area of a triangle made up of edges **a, b, d.** That is, $16 \cdot S_2^2=g^4-2(a^4+b^4+d^4)$. Thus it turns out that $\mathbf{g^4=16 \cdot S_1^2-16 \cdot S_2^2}$.

From this relation it is clear that if all the elements of the cuboid are integers, then the spatial diagonal is even, but according to the conditions of the problem it must be odd. Therefore, one of the elements of the cuboid is a non-integer.

But there's a problem. If the existence of a triangle based on the lateral diagonals (c, f, e) is a fact that does not require proof, then the existence of a triangle based on the edges of a cuboid (a, b, d) requires justification. Examples show that the radical expression in Heron's formula can be both negative and positive, that is, a triangle may or may not exist, but surprisingly, Heron's formula works in both cases.

Using the following example, we will show how many cuboids can be constructed on one spatial diagonal by multiplying two triangles $9^2+40^2=41^2$ и $8^2+15^2=17^2$.

**Example 1:**
1. ($a=104$, $d=153$, $b=672$, $c=680$, $e=185$, $f^2=474993$, $\mathbf{g=697}$).
   Squared area of the triangle (**c, f, e**) $=16 \cdot S_1^2=62834616576$,
   Squared area of the triangle (**a, b, d**) $\to 16 \cdot S_2^2= - 173175767905 \to$
   $g^4 = 16 \cdot S_1^2 - 16 \cdot S_2^2 \to (62834616576) - (-173175767905) = \mathbf{697^4}$.
2. ($a=328$, $b^2=171200$, $d=455$, $e=615$, $c=528$, $f^2=314609$, $\mathbf{g=697}$)
   $16 \cdot S_1^2=304534553600$
   $16 \cdot S_2^2=68524169119 \to g^4 = 16 \cdot S_1^2-16 \cdot S_2^2=304534553600-68524169119=\mathbf{697^4}$.
3. ($a^2=344000$, $b=672$, $d=615$, $c^2=107584$, $e=185$, $f^2=829809$, $\mathbf{g=697}$)
   $16 \cdot S_1^2= (-458616750400)$
   $16 \cdot S_2^2=(-694626134881) \to g^4 =16 \cdot S_1^2-16 \cdot S_2^2=(-458616750400)-(-694626134881)= \mathbf{697^4}$.
4. ($a^2= (i^2 183616)$, $b=680$, $d=455$, $c=528$, $e=153$, $f^2=669425$, $\mathbf{g=697}$).
   $16 \cdot S_1^2= (-108755123200)$
   $16 \cdot S_2^2= (-344765507681) \to g^4= 16 \cdot S_1^2-16 \cdot S_2^2=(-108755123200)-(-344765507681)=\mathbf{697^4}$
5. ($a=320$, $b=600$, $d=135$, $c=680$, $e^2=120625$, $\mathbf{g^2=480625}$)
   $16 \cdot S_1^2=181164960000$
   $16 \cdot S_2^2= (-49835430625) \to 16 \cdot S_1^2-16 \cdot S_2^2=181164960000-(-49835430625) =480625^2$

And so on, in a similar way it is possible to calculate the difference of the squared areas for the subsequent cuboids.

6. (edges 600, 320, 72) lateral diagonals 328, ($\sqrt{365184}$), 680 spatial diagonal ($\sqrt{467584}$)
7. (edges 72, 135, 320) lateral diagonals 328, 153, ($\sqrt{120625}$) spatial diagonal ($\sqrt{125809}$)
8. (edges 72, 135, 600) lateral diagonals 615, 153, ($\sqrt{365184}$) spatial diagonal ($\sqrt{383409}$).

**Example 2:** In this example we present only two cuboids with one non-integer element.
1. ($a=520$, $d^2=618849$, $b=576$, $c=776$, $e=943$, $f=975$, $\mathbf{g=1105}$).
   Squared area of the triangle (**c, f, e**) $16 \cdot S_1^2=1849472983296$.
   Squared area of the triangle (**a, b, d**) $\to 16 \cdot S_2^2=358570932671 \to$
   $\mathbf{g^4} = 16 \cdot S_1^2-16 \cdot S_2^2 \to (1849472983296) - (358570932671) =\mathbf{1105^4}$.
2. ($a=448$, $d=975$, $b=264$, $c=520$, $e=1073$, $f^2=1020321$, $\mathbf{g=1105}$).
   Squared area of the triangle (**c, f, e**) $16 \cdot S_1^2=1084149063936$
   Squared area of the triangle (**a, b, d**) $\to 16 \cdot S_2^2= -406752986689$
   $\mathbf{g^4} = 16 \cdot S_1^2-16 \cdot S_2^2 \to 1084149063938 - (-406752986689) =\mathbf{1105^4}$.

**Trigonometric relations between elements of a rational cuboid**

**Proposition 6.** Let us assume the following:

for $a^2 +b^2 =c^2$ let $\mathrm{ch}^2\, \alpha=\frac{b^2}{b^2 - a^2}$, $\mathrm{sh}^2\alpha=\frac{a^2}{b^2 - a^2}$

for $a^2 +d^2 =e^2$ let $\mathrm{ch}^2\beta=\frac{d^2}{d^2 - a^2}$, $\mathrm{sh}^2\beta=\frac{a^2}{d^2 - a^2}$

for $b^2 + d^2 = f^2$ let $\mathrm{ch}^2\gamma=\frac{b^2}{b^2 - d^2}$, $\mathrm{sh}^2\gamma=\frac{d^2}{b^2 - d^2}$

Then $\mathrm{th}^2\alpha \cdot \mathrm{cth}^2\beta \cdot \mathrm{cth}^2\gamma =1$.

Proof: Instead of functions we substitute their numerical values $(a^2/b^2) \cdot (d^2/a^2) \cdot (b^2/d^2)=1$.

**Lemma 4**. The angles of the triangles (Proposition 6) describing a rational cuboid are related by the following relation: $e^{2\alpha+2\beta}+e^{2\alpha+2\gamma}=e^{2\beta+2\gamma}+1$.

Proof: Let us express $\mathrm{th}\alpha \cdot \mathrm{cth}\beta \cdot \mathrm{cth}\gamma =1$ in terms of hyperbolic functions $e^x$ and $e^{-x}$, after some simple transformations we will arrive at the desired formula.

**Lemma 5.** The angles of the triangles (Proposition 6) describing a rational cuboid can be found using the following formulas:

$$2\alpha= \ln \left[\frac{\mathrm{ch}(\beta+\gamma)}{\mathrm{ch}(\beta - \gamma)}\right], 2\beta= \ln \left[\frac{\mathrm{sh}(\alpha+\gamma)}{\mathrm{sh}(\alpha - \gamma)}\right], 2\gamma= \ln \left[\frac{\mathrm{sh}(\alpha+\beta)}{\mathrm{sh}(\alpha - \beta)}\right].$$

Proof:

$\mathrm{th}^2\alpha \cdot \mathrm{cth}^2\beta \cdot \mathrm{cth}^2\gamma = 1 \to \mathrm{sh}\alpha \cdot \mathrm{ch}\beta \cdot \mathrm{ch}\gamma=\mathrm{ch}\alpha \cdot \mathrm{sh}\beta \cdot \mathrm{sh}\gamma \to$

$$\text{sh}\alpha\cdot\frac{[\text{ch}(\beta+\gamma)+\text{ch}(\beta-\gamma)]}{2} = \text{ch}\alpha\cdot\frac{[\text{ch}(\beta+\gamma)-\text{ch}(\beta-\gamma)]}{2}. \rightarrow$$

$\text{ch}(\beta+\gamma)(\text{sh}\alpha-\text{ch}\alpha)+\text{ch}(\beta-\gamma)(\text{sh}\alpha+\text{ch}\alpha)=0.$ $[(\text{sh}\alpha-\text{ch}\alpha)=-e^{-\alpha}.]\rightarrow$

$$\frac{[-\text{ch}(\beta+\gamma)]}{e}{}^{-\alpha} +\text{ch}(\beta-\gamma)\cdot e^{\alpha} =0. \rightarrow 2\alpha= \ln\left[\frac{\text{ch}(\beta+\gamma)}{\text{ch}(\beta-\gamma)}\right].$$

**Proposition 7.** Let us assume the following:

for $c^2 + d^2 = g^2$ let $\cos^2\theta=\frac{c^2}{c^2+d^2}$, $\sin^2\theta=\frac{d^2}{c^2+d^2}$

for $b^2 + e^2 = g^2$ let $\cos^2\varepsilon=\frac{b^2}{b^2+e^2}$, $\sin^2\varepsilon=\frac{e^2}{b^2+e^2}$

for $a^2 + f^2 = g^2$ let $\cos^2\delta=\frac{f^2}{f^2+a^2}$, $\sin^2\delta=\frac{a^2}{f^2+a^2}$

Then $\cos^2\theta+\sin^2\varepsilon+\cos^2\delta=2$, a $\sin^2\theta+\cos^2\varepsilon+\sin^2\delta=1$.

Proof: Same as in **Proposition 6** we substitute their numerical values.